\documentclass[reqno,12pt]{amsart}
\usepackage{amsfonts,amsmath,amsthm,amssymb,amscd}
\usepackage[all]{xy}
\newcommand{\version}{Ver.~0.0}
\newcommand{\setversion}[1]{\renewcommand{\version}{Ver.~{#1}}}
\setversion{0.0 [2009/06/14 13:00]}
\setversion{0.1 [2009/08/09 18:43]}
\setversion{0.5 [2009/11/11 10:29]}
\setversion{0.6 [2009/11/20 22:35]}
\setversion{0.7 [2009/11/27]}
\setversion{0.8 [2010/02/13 12:57]}
\setversion{0.9 [2010/02/22]}
\setversion{1.0 [2010/03/02 07:42]}
\setversion{1.1 [2010/03/07 05:04]}
\setversion{1.2 [2010/03/07 18:11]}
\setversion{1.3 [2010/04/01 16:56]}
\setversion{2.1 [2010/08/31]}

\title [Double flag variety]
{Double flag varieties for a symmetric pair and finiteness of orbits}

\author{Kyo Nishiyama}
\address{
Department of Physics and Mathematics\\
Aoyama Gakuin University\\
Fuchinobe 5-10-1, Sagamihara 229-8558, Japan}
\email{kyo@gem.aoyama.ac.jp}

\thanks{Supported by JSPS Grant-in-Aid for Scientific Research (B) \#{21340006}.}

\author{Hiroyuki Ochiai}
\address{
Faculty of Mathematics\\
Kyushu University\\
744, Motooka, Nishi-ku, Fukuoka 819-0395, Japan}
\email{ochiai@math.kyushu-u.ac.jp}

\thanks{Supported by JSPS Grant-in-Aid for Scientific Research (A) \#{19204011}.}

\date{}

\theoremstyle{plain}
\newtheorem{theorem}{Theorem}
\newtheorem{proposition}[theorem]{Proposition}
\newtheorem{corollary}[theorem]{Corollary}
\newtheorem{lemma}[theorem]{Lemma}

\newtheorem{problem}[theorem]{Problem}
\newtheorem{itheorem}{Theorem}

\theoremstyle{definition}

\theoremstyle{remark}

\numberwithin{equation}{section}
\numberwithin{theorem}{section}

\newcommand{\R}{\mathbb{R}}

\newcommand{\bbK}{\mathbb{K}}

\newcommand{\bbG}{\mathbb{G}}
\newcommand{\C}{\mathbb{C}}
\newcommand{\bbP}{\mathbb{P}}
\newcommand{\bbQ}{\mathbb{Q}}

\newcommand{\lie}[1]{\mathfrak{#1}}
\newcommand{\lier}[1]{\mathfrak{#1}_{\R}}

\newcounter{thmenum}
\newenvironment{thmenumerate}{%
\begin{list}{$(\thethmenum)$}{%
\usecounter{thmenum}
\setlength{\labelsep}{.5em}
\setlength{\labelwidth}{-7pt}
\setlength{\topsep}{0pt}
\setlength{\partopsep}{0pt}
\setlength{\parsep}{0pt}
\setlength{\leftmargin}{3pt}
\setlength{\rightmargin}{0pt}
\setlength{\itemindent}{\leftmargin}
\setlength{\itemsep}{0pt}
}}
{\end{list}}

\newcommand{\mycomment}[1]{} 

\newlength{\lengthcup}
\newcommand{\Hom}{\qopname\relax o{Hom}}

\newcommand{\Image}{\qopname\relax o{Im}}

\newcommand{\closure}[1]{\overline{#1}}

\newcommand{\restrict}{\big|}

\newcommand{\orbit}{\mathbb{O}}
\newcommand{\calorbit}{\mathcal{O}}
\newcommand{\irreps}[1]{\mathrm{Irr}(#1)}

\newcommand{\GL}{\mathrm{GL}}
\newcommand{\SL}{\mathrm{SL}}

\newcommand{\SO}{\mathrm{SO}}
\newcommand{\Sp}{\mathrm{Sp}}

\newcommand{\Grass}{\qopname\relax o{Grass}}
\newcommand{\LGrass}{\qopname\relax o{LGrass}}
\newcommand{\IGrass}{\qopname\relax o{IGrass}}

\newcommand{\injection}{\hookrightarrow}

\newcommand{\GFl}{\mathfrak{X}}
\newcommand{\affGFl}{\widehat{\mathfrak{X}}}
\newcommand{\GFltheta}{\mathfrak{X}^{\theta}}
\newcommand{\KFl}{\mathcal{Z}}
\newcommand{\GKFl}[2]{\GFl_{#1}\times\KFl_{#2}}

\newcommand{\Aqm}[2]{\mathcal{A}_{#1}({#2})}

\newcommand{\lengthof}[1]{\ell(#1)}

\newcommand{\SpD}{{\Sp}D}
\newcommand{\SpE}{{\Sp}E}
\newcommand{\SpY}{{\Sp}Y}

\newcommand{\opp}[1]{#1^{\,\circ}}

\newcommand{\inclusion}{\hookrightarrow}

\newcommand{\Dtheta}{\Delta_{\theta}}

\newcommand{\thetacalorbit}{\calorbit^{\theta}}

\newcommand{\llg}{\lie{g}}
\newcommand{\llk}{\lie{k}}
\newcommand{\llp}{\lie{p}}
\newcommand{\lls}{\lie{s}}

\newcommand{\miniP}{P_{\mathrm{min}}}
\newcommand{\maxT}{\mathcal{T}}
\newcommand{\maxTtheta}{\maxT^{\theta}}
\newcommand{\CTX}{\mathcal{C}}
\newcommand{\CthetaTX}{\mathcal{C}^{\theta}}
\newcommand{\RSV}{\mathcal{V}}
\newcommand{\TwistedInv}{\mathcal{I}}

\newcommand{\bsl}{\backslash}

\begin{document}

\begin{abstract}
Let $ G $ be a reductive algebraic group over the complex number filed, 
and $ K = G^{\theta} $ be the fixed points of an involutive automorphism $ \theta $ of $ G $ 
so that $ (G, K) $ is a symmetric pair.   

We take parabolic subgroups $ P $ and $ Q $ of $ G $ and $ K $ respectively and  
consider a product of partial flag varieties $ G/P $ and 
$ K/Q $ with diagonal $ K $-action.  
The double flag variety $ G/P \times K/Q $ thus obtained 
is said to be \emph{of finite type} if there are finitely many $ K $-orbits on it.  
A triple flag variety $ G/P^1 \times G/P^2 \times G/P^3 $ is a special case of 
our double flag varieties, and there are many interesting works on the triple flag varieties.

In this paper, we study double flag varieties $ G/P \times K/Q $ of finite type.  
We give efficient criterion under which the double flag variety is of finite type.  
The finiteness of orbits is strongly related to spherical actions of $ G $ or $ K $.  
For example, we show a partial flag variety $ G/P $ is $ K $-spherical 
if a product of partial flag varieties $ G/P \times G/\theta(P) $ is $ G $-spherical.
We also give many examples of the double flag varieties of finite type, and for type AIII, 
we give a classification when $ P = B $ is a Borel subgroup of $ G $. 
\end{abstract}

\maketitle


\section*{Introduction}

Let $ G $ be a connected reductive algebraic group over $ \C $.  
Recently, there appear many interesting examples of product of 
(partial) flag varieties which have finitely many $G$-orbits.  One 
example is $ X = G/B \times G/B \times \mathbb{P}(V) $ where $ G = 
GL(V) $, $ B $ a Borel subgroup and $ \mathbb{P}(V) $ denotes the 
projective space over $ V $.  The third factor $ \mathbb{P}(V)$ is isomorphic to a partial flag variety 
$ G/P $, where $ P $ is a maximal parabolic subgroup stabilizing a one dimensional subspace of $ V $.   
It is known that there are finitely many $ G $-orbits on 
$ X $, and by the work of Travkin, Finkelberg and Ginzburg (\cite{Travkin.2009, FGT.2009}), 
there are miraculous similarities between $ X $ and Steinberg variety.  For 
example, they established a kind of Robinson-Schensted-Knuth 
correspondence for the orbits on $ X $, and study some micro-local 
properties using Hecke algebras.  
The maximal parabolic $ P $ above is called ``mirabolic" after Ginzburg.  

In general, one can consider a triple product of partial flag varieties.  
For a parabolic subgroup $ P $ of $ G $, 
we denote $ \GFl_{P} = G/P $ a partial flag variety.  
Magyar-Weymann-Zelevinsky (\cite{MWZ.1999, MWZ.2000}) classified 
the triple flag varieties $ \GFl_{P^1} \times \GFl_{P^2} \times \GFl_{P^3} $ 
which have finitely many $ G $-orbits when $ G $ is a classical group of type A or type C.  
They also gave parametrizations of orbits.

In this paper, we generalize the notion of triple flag varieties to a symmetric pair 
$ (G, K) $, where $ K $ is a symmetric subgroup of $ G $ consisting of the fixed points 
of an involutive automorphism $ \theta $.  
Thus we take parabolic subgroups $ P \subset G $ and $ Q \subset K $, 
and consider a product of partial flag varieties $ \GFl_{P} = G/P $ and $ \KFl_{Q} = K/Q $.  
The group $ K $ acts on the \emph{double flag variety} $ \GFl_{P} \times \KFl_{Q} $ diagonally.  

If one considers $ \bbG = G \times G $ and an involution $ \theta(g_1, g_2) = (g_2, g_1) $ of $ \bbG $, 
the symmetric subgroup $ \bbK = \bbG^{\theta} $ is just the diagonal subgroup 
$ \Delta(G) \subset \bbG $.  
Then, for parabolic subgroups $ \bbP = (P^1, P^2) \subset \bbG $ and $ \bbQ = \Delta(P^3) \subset \bbK $, 
the double flag variety can be interpreted as 
\begin{equation*}
\bbG/ \bbP \times \bbK/ \bbQ
= (G \times G)/(P^1 \times P^2) \times \Delta(G)/\Delta(P^3)
\simeq \GFl_{P^1} \times \GFl_{P^2} \times \GFl_{P^3}
\end{equation*}
which is nothing but the triple flag variety.  
So our double flag variety is a generalization of triple flag varieties.

We say a double flag variety $ \GFl_{P} \times \KFl_{Q} $ is of \emph{finite type} if 
there are only finitely many $ K $-orbits on it.  
One of the most interesting problems is to classify the double flag varieties of finite type.  
We give two efficient criterions for the finiteness of orbits using triple flag varieties.  
Namely, in Theorem 
\ref{theorem:finiteness.of.3.flags.implies.finiteness.of.2.flags}, we prove 

\begin{itheorem}
\label{itheorem:triple.flag.to.double.flag}
Let $ P' $ be a $ \theta $-stable parabolic of $ G $ such that $ P' \cap K = Q $.  
If the number of $ G $-orbits on $ \GFl_{P} \times \GFl_{\theta(P)} \times \GFl_{P'} $ is finite, 
then there are only finitely many $ K $-orbits on the double flag variety $ \GKFl{P}{Q} $.
\end{itheorem}

The next theorem (Theorem \ref{proposition:finiteness.by.HC.embedding}) is also useful.

\begin{itheorem}
\label{iproposition:intersection.of.parabolics}\label{itheorem:intersection.of.parabolics}
Let $ P^i \; (i = 1, 2, 3) $ be a parabolic subgroup of $ G $.  
Suppose that $\GFl_{P^1} \times \GFl_{P^2} \times \GFl_{P^3}$
has finitely many $G$-orbits
and that $Q:= P^2 \cap P^3$ is a parabolic subgroup of $K$.
Then $\GKFl{P^1}{Q}$ has finitely many $K$-orbits.  

Moreover, if $ P^1 $ is a Borel subgroup $ B $ and the product $ P_2 P_3 $ is open in $ G $, 
then the converse is also true, i.e., 
the double flag variety 
$ \GKFl{B}{Q} $ 
is of finite type 
if and only if 
the triple flag variety 
$ \GFl_{B} \times \GFl_{P^2} \times \GFl_{P^3} $ 
is of finite type.
\end{itheorem}

We construct many examples of double flag varieties of finite type using 
Theorems \ref{itheorem:triple.flag.to.double.flag} and 
\ref{iproposition:intersection.of.parabolics} in \S\S~\ref{section:typeA}--\ref{section:typeC}, 
and if $ P = B $ is Borel, we give complete classification for certain cases.
However, in general cases, the classification of double flag varieties of finite type 
seems to be difficult.  

Double flag varieties of finite type are strongly related to {spherical actions} of 
$ G $ or $ K $.  
Recall that an action of a reductive algebraic group $ G $ on a variety $ X $ is called \emph{spherical} 
if there is an open dense $ B $-orbit for a certain Borel subgroup $ B $ of $ G $.  
The existence of an open dense $ B $-orbit is in fact equivalent to 
the finiteness of $ B $-orbits on $ X $ 
due to Brion \cite[\S~1.5]{Brion.PM80.1989} and independently to 
Vinberg \cite{Vinberg.1986}.  
We often use this finiteness criterion for spherical actions below.  

The following theorem (Theorem \ref{theorem:Kspherical.partial.flags}), 
which is a corollary of the first theorem, 
exhibits a good connection to the spherical action.

\begin{itheorem}
\label{itheorem:K.spherical.PFV}
Let $ P $ be a parabolic subgroup of $ G $.  
If $ \GFl_P \times \GFl_{\theta(P)} $ is a spherical $ G $-variety, then 
$ \GFl_P $ is a spherical $ K $-variety.
\end{itheorem}

For a parabolic subgroup $ P $ in $ G $, we can find 
a finite-dimensional irreducible representation  $ V_{\lambda} $ of $ G $ with highest weight $ \lambda $ 
such that $ P = \{ g \in G \mid g v_{\lambda} \in \C v_{\lambda} \} $, where $ v_{\lambda} $ denotes 
a highest weight vector of $ V_{\lambda} $.  
Assume that the conclusion of Theorem \ref{itheorem:K.spherical.PFV} holds, i.e., 
we assume that $ \GFl_P $ is $ K $-spherical.  
Then the contragredient $ V_{k \lambda}^{\ast}\restrict_K $ decomposes multiplicity freely as a $ K $-module for any 
non-negative integer $ k \geq 0 $ (see Lemma \ref{lemma:K.spherical.PFV.and.K.MF}).  
This is one of interesting conclusions of Theorem \ref{itheorem:K.spherical.PFV}.

We will also give some other examples of spherical multiple flag varieties in \S~\ref{section:spherical.action}.  

There seems to be intimate connection between double flag varieties of finite type and visible actions 
(see \cite{Kobayashi.2005} for the definition of visible actions).  
Let us denote the compact real form of $ K $ by $ K_U $.  
Then $ K_U $ acts on $ \GFl_{P} = G/P $ visibly if and only if $ \GFl_{P} $ is $ K $-spherical.  
This is equivalent to say that $ \GKFl{P}{S} $ is of finite type for a Borel subgroup $ S $ of $ K $.  
See also \cite[Cor.~17]{Kobayashi.2005} and \cite{Kobayashi.2008}.

\smallskip
\noindent{\it Acknowledgment.} 
We started the study of double flag varieties when we visited NCTS in National Cheng-Kung University in Tainan.  
We thank Ngau Lam for his generous hospitality at NCTS.  
We also thank Peter Trapa for useful discussion.  

Finally, we thank the anonymous referee, who carefully read the earlier version of the paper.  
The references \cite{Vinberg.1986} and \cite{Kramer.1976} are due to him/her.  
Also referee's opinions have considerably improved the whole organization of the paper.

\section{Double flag varieties for symmetric pair}

Let $ G $ be a connected reductive algebraic group over the complex number field $ \C $, and $ \theta $ its (non-trivial) involutive automorphism.  
We put $ K = G^{\theta} $, a subgroup whose elements are fixed by $ \theta $, 
and call it a symmetric subgroup of $ G $.  
We denote the Lie algebra of $ G $ (respectively of $ K $) by $ \lie{g} $ (respectively $ \lie{k} $).  
In the following, we use the similar notation; 
for an algebraic group we use a Roman capital letter, and for its Lie algebra the corresponding German small letter.

For a parabolic subgroup $ P $, we denote a partial flag variety 
consisting of all $ G $-conjugates of $ P $ by $ \GFl_P $.  
Since $ P $ is self-normalizing, $ \GFl_P $ is isomorphic to $ G/P $ as a $ G $-variety.  
We also choose a $ \theta $-stable parabolic $ P' $ in $ G $, and put $ Q = K \cap P' $.  
Then $ Q $ is a parabolic subgroup of $ K $, and every parabolic subgroup of $ K $ can be obtained in this way 
(see \cite[Theorem~2]{Brion.Helminck.2000}).
We denote a partial flag variety $ K/Q $ by $ \KFl_Q $.

We consider the following problem.

\begin{problem}
\label{problem:finite.K.orbits.on.double.flag}
Let the symmetric subgroup $ K $ act on the product of the partial flag varieties $ \GKFl{P}{Q} $ diagonally. 
\begin{thmenumerate}
\item
\label{prob:i}
Classify all the pair $ (P, Q) $ {\upshape(}or $ (P, P') ${\upshape)} for a given pair $ (G, K) $ 
which admits finitely many $ K $-orbits on $ \GKFl{P}{Q} $.  
We are also interested in the case where $ \GKFl{P}{Q} $ contains an open $ K $-orbit.  
\item
\label{prob:ii}
If there are finitely many orbits, classify all the $ K $-orbits on $ \GKFl{P}{Q} $, and 
study the geometry of orbits; for example, closure relations, combinatorial descriptions, equivariant cohomology and so on.
\item
\label{prob:iii}
Establish a relation to the representation theory of Harish-Chandra $ ( \lie{g}, K) $-modules.
\end{thmenumerate}
\end{problem}

All these problems are still open, but at the same time many (explicit) results are already known.  
We will give partial answers to the problem \eqref{prob:i}, which are new, in the following sections.

The problem \eqref{prob:iii} may require some account.  Let us explain it briefly.  
Since we assume $ P' $ is a $ \theta $-stable parabolic subgroup and $ Q = K \cap P' $, 
$ \KFl_Q = K/Q $ can be embedded into $ \GFl_{P'} = G/P' $, which is a partial flag variety of $ G $.  
In fact, $ \KFl_Q $ is isomorphic to a closed $ K $-orbit in $ \GFl_{P'} $, and for every closed $ K $-orbit in $ \GFl_{P'} $, 
one can attach a Harish-Chandra $ ( \lie{g}, K ) $-module via Beilinson-Bernstein theory.  
This module is known to be a derived functor module $ \Aqm{\lie{p}'}{\rho'} $ induced 
from a one-dimensional character of a $ \theta $-stable 
parabolic subalgebra $ \lie{p}' $ contained in the $ K $-orbit.  
On the other hand, there is an open dense $ K $-orbit on $ \GFl_P $ which should correspond to a degenerate principal series representation if some translate of $ \lie{p} $ has a real form.  
Thus, in a suitable space of cohomology of an invertible sheaf over $ \GKFl{P}{Q} $, 
one can hopefully realize the tensor product of a degenerate principal series representation and 
$ \Aqm{\lie{p}'}{\rho'} $.
Our condition of the finiteness of the $ K $-orbits put a restriction 
on the tensor product and we expect a kind of multiplicity-free property on the tensor product.

\section{Triple flags}

Let us return to Problem \ref{problem:finite.K.orbits.on.double.flag} and 
consider a symmetric pair 
$ (\bbG, \bbK) $, where $ \bbG = G \times G $ for a reductive group $ G $ over $ \C $ and 
$ \bbK = \Delta G $ is the diagonal embedding.  
This symmetric subgroup $ \bbK $ corresponds to the involution 
$ \theta : \bbG \to \bbG $ defined by $ \theta(g_1, g_2) = (g_2, g_1) $.  
Take a parabolic subgroup $ \bbP = P^1 \times P^2 $ in $ \bbG $, where $ P^i $ is a parabolic subgroup of $ G $.  
A $ \theta $-stable parabolic subgroup $ \bbP' $ in $ \bbG $ can be written as 
$ \bbP' = P^3 \times P^3 $ for a certain parabolic subgroup $ P^3 $ in $ G $.  
Thus $ \bbQ = \bbK \cap \bbP' = \Delta P^3 $ is a diagonal subgroup in $ \Delta G $.  
Now it is immediate to see that, in this setting, our Problem \ref{problem:finite.K.orbits.on.double.flag} can be translated into 

\begin{problem}
\label{problem:finite.G.orbits.on.triple.flag}
Let $ G $ act on the triple product of partial flag varieties $ \GFl_{P^1} \times \GFl_{P^2} \times \GFl_{P^3} $ diagonally. 
\begin{thmenumerate}
\item
\label{prob.G:i}
Classify all the triples $ (P^1, P^2, P^3) $, for which there are finitely many $ G $-orbits on the triple product $ \GFl_{P^1} \times \GFl_{P^2} \times \GFl_{P^3} $.  
If this is the case, we say the triple product is \emph{of finite type}.
\item
\label{prob.G:ii}
If there are finite number of orbits, classify all the $ G $-orbits and 
study the geometry of orbits.
\item
\label{prob.G:iii}
Establish a relation to the representation theory.
\end{thmenumerate}
\end{problem}

This problem (at least for \eqref{prob.G:i} and \eqref{prob.G:ii}) was solved almost completely for classical groups 
by Magyar-Weymann-Zelevinsky (\cite{MWZ.1999, MWZ.2000}), Travkin (\cite{Travkin.2009}) and Finkelberg-Ginzburg-Travkin (\cite{FGT.2009}).
If $ P^3 $ is a Borel subgroup, Littelmann (\cite{Littelmann.1994}) investigated a representation theoretic meaning; in fact, 
it seems that this work is one of motivations of \cite{MWZ.1999, MWZ.2000}.  
We do not go into the details of their works, 
but let us introduce the classification achieved by \cite{MWZ.1999, MWZ.2000} without proof since we need it later.

\subsection{Type A}
\label{subsection:MWZ.typeA}

Let $ G = \GL_n(\C) $ be the general linear group, 
which we will denote simply by $ \GL_n $ if there is no confusion.  
To specify a parabolic subgroup $ P $ of $ G $, we use an unordered partition (or composition) 
 $ \lambda $ of $ n $; 
i.e., if $ P = P_{\lambda} $ corresponds to $ \lambda = (\lambda_1, \lambda_2, \dots, \lambda_l) $, its Levi part is in block diagonal form of 
$ \GL_{\lambda_1} \times \GL_{\lambda_2} \times \cdots \times \GL_{\lambda_l} $, and its unipotent radical is in upper triangular form.  
The number of non-zero parts in $ \lambda $ is denoted by $ \lengthof{\lambda} $ and is called the length of $ \lambda $.

\begin{theorem}[Magyar-Weymann-Zelevinsky]
\label{theorem:MWZtypeA}
Let $ \GFl_{P} = G/P $ be a partial flag variety, where $ G = \GL_n $.  
\begin{thmenumerate}
\item
For a collection of proper parabolic subgroups $ P^1, \dots, P^k $, 
if the number of $ G $-orbits on $ \GFl_{P^1} \times \GFl_{P^2} \times \cdots \times \GFl_{P^k} $ is finite, then $ k \leq 3 $.
\item
A triple product $ \GFl_{P_{\lambda}} \times \GFl_{P_{\mu}} \times \GFl_{P_{\nu}} $ of partial flag varieties
is of finite type if and only if it is from the following list 
(with possible change of the order of parabolic subgroups, and the order of parts of partitions involved).
\begin{equation*}
\begin{array}{ccccc}
\text{type} & (\lengthof{\lambda}, \lengthof{\mu}, \lengthof{\nu}) & \text{extra condition(s)} \\ \hline
S_{q, r} & (2, q, r) & \lambda = (n - 1, 1) \\
D_{r + 2} & (2, 2, r) & \\
E_6 & ( 2, 3, 3) & \\
E_7 & ( 2, 3, 4) & \\
E_8 & ( 2, 3, 5) & \\
E_{r + 3}^{(a)} & ( 2, 3, r) & \lambda = (n - 2, 2) \; (n \geq 4) \\
E_{r + 3}^{(b)} & ( 2, 3, r) & \mu = (\mu_1, \mu_2, 1) 
\end{array}
\end{equation*}
\end{thmenumerate}
\end{theorem}

For the first statement, note that if $ k = 1 $ then $ \GFl_P $ is homogeneous; 
if $ k = 2 $, then $ G \backslash (\GFl_{P^1} \times \GFl_{P^2}) \simeq P^1 \backslash G / P^2 $, which is 
further isomorphic to $ W_{P^1} \backslash W / W_{P^2} $ by the Bruhat decomposition (we denote by $ W $ the Weyl group of $ G $, and 
by $ W_{P} $ the Weyl group of $ P $).  So they are always of finite type.

\subsection{Type C}
\label{subsection:MWZ.typeC}

In this subsection we put $ G = \Sp_{2n}(\C) $, which we abbreviate to $ \Sp_{2n} $.  
The symplectic group $ G $ acts on partial flags of isotropic subspaces 
$ F_1 \subset F_2 \subset \dots \subset F_l $ of fixed dimensions.  
Let us denote the orthogonal subspace of $ F_i $ by $ F_i^{\bot} $.  
Then we have a partial flag of subspaces
\begin{equation}
\label{eqn:flag.of.typeC}
F_1 \subset F_2 \subset \cdots \subset F_{l - 1} \subset F_l 
\subset F_l^{\bot} \subset F_{l - 1}^{\bot} \subset \cdots \subset F_2^{\bot} \subset F_1^{\bot} .
\end{equation}
A parabolic subgroup of $ G $ is specified as a fixed point subgroup of a partial flag as in \eqref{eqn:flag.of.typeC}, 
and its conjugacy class is determined by the dimensions of the subspaces in the flag.  
We put, if $ \dim F_l < n $, then 
\begin{equation*}
\begin{cases}
\lambda_i = \dim F_i/ F_{i - 1} & (1 \leq i \leq l) \\
\lambda_{l + 1} = \dim F_l^{\bot}/ F_l = 2 (n - \dim F_l ) & \\
\lambda_{l + i + 1} = \dim F_{l - i}^{\bot}/ F_{l - i + 1}^{\bot} = \lambda_{l - i + 1} & (1 \leq i \leq l) 
\end{cases}
\end{equation*}
with $ F_0 $ understood as $ \{ 0 \} $,  
and if $ \dim F_l = n $ 
\begin{equation*}
\begin{cases}
\lambda_i = \dim F_i/ F_{i - 1} & (1 \leq i \leq l) \\
\lambda_{l + i} = \dim F_{l - i}^{\bot}/ F_{l - i + 1}^{\bot} = \lambda_{l - i + 1} & (1 \leq i \leq l) 
\end{cases}
\end{equation*}
Then $ \lambda = ( \lambda_1, \lambda_2, \dots ) $ is an unordered partition with $ |\lambda| = 2 n $, 
where $ |\lambda| $ is the size of $ \lambda $.  
We denote by $ P = P_{\lambda} $ the corresponding parabolic subgroup, whose Levi part is isomorphic to 
$ \GL_{\lambda_1} \times \cdots \times \GL_{\lambda_l} \times \Sp_{\lambda_{l + 1}} $.  
Here the factor $ \Sp_{\lambda_{l + 1}} $ does not appear if $ \dim F_l = n $.  
Thus, if $ \lambda = (n, n) $, the corresponding parabolic subgroup $ P_{(n, n)} $ is a Siegel parabolic with Levi component $ \GL_n $; and 
if $ \lambda = (m, 2 (n - m), m) $ with $ m < n $, then $ P_{(m,  2(n - m), m)} $ is a maximal parabolic subgroup with Levi component $ \GL_m \times \Sp_{2(n - m)} $.

With this notation, we can state the following

\begin{theorem}[Magyar-Weymann-Zelevinsky]
\label{theorem:MWZtypeC}
Let $ \GFl_{P} = G/P $ be a partial flag variety, where $ G = \Sp_{2n} $.  
\begin{thmenumerate}
\item
For a collection of proper parabolic subgroups $ P^1, \dots, P^k $, 
if the number of $ G $-orbits on $ \GFl_{P^1} \times \GFl_{P^2} \times \cdots \times \GFl_{P^k} $ is finite, then $ k \leq 3 $.
\item
A triple product $ \GFl_{P_{\lambda}} \times \GFl_{P_{\mu}} \times \GFl_{P_{\nu}} $ of partial flag varieties 
is of finite type if and only if it is from the following list 
(up to appropriate changes of the order of parabolic subgroups, and the order of parts of partitions involved).
\begin{equation*}
\begin{array}{ccccc}
\text{type} & (\lengthof{\lambda}, \lengthof{\mu}, \lengthof{\nu}) & \text{extra condition(s)} \\ \hline
\SpD_{r + 2} & (2, 2, r) & \lambda = \mu = (n, n) \\
\SpE_6 & ( 2, 3, 3) & \lambda = (n, n) \\
\SpE_7 & ( 2, 3, 4) & \lambda = (n, n) \\
\SpE_8 & ( 2, 3, 5) & \lambda = (n, n) \\
\SpE_{r + 3}^{(b)} & ( 2, 3, r) & \lambda = (n, n); \; \mu = (1, 2 n - 2, 1); \; 3 \leq r \\
\SpY_{4, r} & (3, 3, r) & \lambda = \mu = (1, 2 n - 2, 1); \; 3 \leq r
\end{array}
\end{equation*}
Here $ \ell(\lambda) $ denotes the length of $ \lambda $.  
So, if $ \ell(\lambda) = 2 $, it implies $ \lambda = (n, n) $ and $ P_{\lambda} $ is a Siegel parabolic.
If $ \ell(\lambda) = 3 $, then $ P_{\lambda} $ is a maximal parabolic as explained above.
\end{thmenumerate}
\end{theorem}

Note that in the case of type C, two of the parabolic subgroups among three should be maximal.  
Moreover, one of those maximal parabolic subgroups must be a Siegel parabolic $ P_{(n, n)} $ 
or $ P_{(1, 2 n - 2 , 1)} $ with Levi component $ \C^{\times} \times \Sp_{2n -2} $.  


\section{Double flag varieties of finite type}
\label{section:double.flags.for.s.pair}

Now we return to our original setting, i.e., 
$ G $ is a reductive group with an involution $ \theta $, and 
$ K = G^{\theta} $ a symmetric subgroup, which is automatically reductive.  
We take a parabolic subgroup $ P $ and a $ \theta $-stable parabolic subgroup $ P' $ in $ G $, 
and put $ Q = P' \cap K $, which is a parabolic subgroup of $ K $.  
As we have already mentioned, 
for any parabolic $ Q \subset K $, we can choose a $ \theta $-stable parabolic subgroup $ P' $ in $ G $ 
which cuts out $ Q $ from $ K $.  
So our assumption causes no essential restriction.

Let us consider the diagonal action of $ K $ on the product of 
partial flag varieties $ \GKFl{P}{Q} $, where $ \GFl_P = G/P $ and $ \KFl_Q = K/Q $.  
We also put $ \GFltheta_{P} = \GFl_{\theta(P)} = G / \theta(P) $.  

The following is one of our main results in this article.

\begin{theorem}
\label{theorem:finiteness.of.3.flags.implies.finiteness.of.2.flags}
If the number of $ G $-orbits on $ \GFl_{P} \times \GFltheta_{P} \times \GFl_{P'} $ is finite, 
then there are only finitely many $ K $-orbits on the double flag variety $ \GKFl{P}{Q} $.
\end{theorem}

\begin{proof}
If $ P' = G $, then this theorem reduces to the well-known fact that 
there are finitely many $ K $-orbits on the partial flag variety $ \GFl_{P} $ 
(\cite{Matsuki.1979, Matsuki.1982}, \cite{Rossmann.1979}, \cite{Springer.1985}).  
For this, there is a beautiful proof by Mili\v{c}i\'{c} \cite[\S~H.2, Theorem 1]{Milicic.1993}, 
and our proof is an extension of his idea to the case of double flags.

Let us consider the following $ \theta $-twisted diagonal embedding:
\begin{equation*}
\Dtheta : \GFl_P \ni P^1 \mapsto ( P^1, \theta(P^1) ) \in \GFl_P \times \GFltheta_P,
\end{equation*}
where we identify $ \GFl_P $ with the set of parabolic subgroups of $ G $ which are conjugate to $ P $.  
Note that $ \theta(P^1) $ belongs to $ \GFltheta_P $ for any $ P^1 \in \GFl_P $.  
Thus we can embed
\begin{equation}
\GKFl{P}{Q} \xrightarrow{\;\sim\;} \Dtheta(\GFl_P) \times \KFl_Q \inclusion \Dtheta(\GFl_P) \times \GFl_{P'} 
\subset \GFl_{P} \times \GFltheta_{P} \times \GFl_{P'} .
\end{equation}
This is a closed embedding, and clearly $ K $-equivariant.
Let us consider a $ \theta $-twisted action of $ G $ on $ \GFl_{P} \times \GFltheta_{P} \times \GFl_{P'} $:
\begin{equation*}
g ( P^1, P^2, {P^3}') = ( g \cdot P^1, \theta(g) \cdot P^2 , g \cdot {P^3}' ) \qquad
(g \in G ; \;\; ( P^1, P^2, {P^3}') \in \GFl_{P} \times \GFltheta_{P} \times \GFl_{P'} ), 
\end{equation*}
which preserves $ \Dtheta(\GFl_P) \times \GFl_{P'} $.  
Note that $ g \in G $ acts on $ \GFl_P $ by conjugation $ g \cdot P^1 = g P^1 g^{-1} $.  
If we indicate this action also by $ \Dtheta $, 
there are only finitely many $ \Dtheta(G) $-orbits on $ \Dtheta(\GFl_P) \times \GFl_{P'} $, namely we have
\begin{equation*}
\Dtheta(G) \backslash \bigl( \Dtheta(\GFl_P) \times \GFl_{P'} \bigr) \simeq 
G \backslash \bigl( \GFl_P \times \GFl_{P'} \bigr) \simeq W_P \backslash W / W_{P'} ,
\end{equation*}
where the last isomorphism comes from the Bruhat decomposition (see \S \ref{subsection:MWZ.typeA}).  
So pick a $ \Dtheta(G) $-orbit $ \thetacalorbit_w $ in $ \Dtheta(\GFl_P) \times \GFl_{P'} $ indexed by 
$ w \in W_P \backslash W / W_{P'} $.  

\begin{lemma}
\label{lemma:finite.orbits.on.O.cap.Otheta}
Take any $ G $-orbit $ \calorbit $ in $ \GFl_{P} \times \GFltheta_{P} \times \GFl_{P'} $, 
and put $ X = \Dtheta(\GFl_P) \times \KFl_Q $.
\begin{thmenumerate}
\item
There are finitely many $ K $-orbits in $ \calorbit \cap \thetacalorbit_w \cap X $.
\item
Let us write the $ K $-orbit decomposition as $ \calorbit \cap \thetacalorbit_w \cap X = \sqcup_{i = 1}^{\ell} \orbit_i $.  
Then each orbit $ \orbit_i $ is a connected component of $ \calorbit \cap \thetacalorbit_w \cap X $, which is 
also an irreducible component as an algebraic variety.
\end{thmenumerate}
\end{lemma}

Let us assume the above lemma.  
Since the decomposition 
\begin{equation*}
X = \bigsqcup_{w \in W_P \backslash W / W_{P'}} \thetacalorbit_w \cap X 
\end{equation*}
is finite, 
and there are only finitely many possibilities of $ G $-orbits $ \calorbit $ in $ \GFl_{P} \times \GFltheta_{P} \times \GFl_{P'} $ 
by the assumption of the theorem, 
we can conclude that $ \# K \backslash X = \# K \backslash ( \GKFl{P}{Q} ) < \infty $.  

Thus it is sufficient to prove the lemma.

Pick a point 
\begin{equation*}
\xi = ( P^1, \theta(P^1), {P^3}') \in \calorbit \cap \thetacalorbit_w \cap X \subset \GFl_{P} \times \GFltheta_{P} \times \GFl_{P'} .  
\end{equation*}
and consider $ \orbit = K \cdot \xi $, a $ K $-orbit through $ \xi $.  
The tangent space of $ \orbit $ at $ \xi $ is contained in 
\begin{equation}
\label{eqn:tangent.of.orbit.subset.intersection}
T_{\xi} \orbit \subset T_{\xi}( \calorbit \cap \thetacalorbit_w \cap X ) \subset T_{\xi} \calorbit \cap T_{\xi} \thetacalorbit_w \cap T_{\xi} X .
\end{equation}
We know 
\begin{align*}
T_{\xi} \calorbit &= \{ ( y + \lie{p}_1, y + \theta(\lie{p}_1), y + \lie{p}_3' ) \mid y \in \lie{g} \} , \\
T_{\xi} \thetacalorbit_w &= \{ ( x + \lie{p}_1, \theta(x) + \theta(\lie{p}_1), x + \lie{p}_3' ) \mid x \in \lie{g} \} .
\end{align*}
We denote a Cartan decomposition by $ \lie{g} = \lie{k} \oplus \lls $, 
where $ \lls $ is the $ (-1) $-eigenspace of the involution $ \theta $ on $ \llg $.    
Let us prove that 
\begin{equation}
\label{eqn:intersection.of.tangents}
T_{\xi} \calorbit \cap T_{\xi} \thetacalorbit_w 
\subset \{ ( z + \lie{p}_1, z + \theta(\lie{p}_1), z + \lls + \lie{p}_3' ) \mid z \in \lie{k} \} .
\end{equation}
To deduce this containment, take a vector from the left hand side.  Then it is expressed as 
\begin{equation*}
( y + \lie{p}_1, y + \theta(\lie{p}_1), y + \lie{p}_3' ) = ( x + \lie{p}_1, \theta(x) + \theta(\lie{p}_1), x + \lie{p}_3' ) 
\end{equation*}
for some $ x \in \llg $ and $ y \in \llg $.  
{From} this, we get
\begin{equation*}
x - y \in \llp_1 , \quad
\theta(x) - y \in \theta(\llp_1) , \quad \text{ and } \quad
x - y \in \llp_3' .
\end{equation*}
By the second formula, we know $ x - \theta(y) \in \llp_1 $, and thus $ y - \theta(y) \in \llp_1 \cap \theta(\llp_1) $.  
Let us decompose $ y $ along the Cartan decomposition:
\begin{equation*}
y = \frac{1}{2} ( y + \theta(y) ) + \frac{1}{2} ( y - \theta(y) ) =: z + v \in \llk \oplus \lls .
\end{equation*}
Then we know $ v \in \llp_1 \cap \theta(\llp_1) $ and 
\begin{equation*}
( y + \lie{p}_1, y + \theta(\lie{p}_1), y + \lie{p}_3' ) 
= ( z + \lie{p}_1, z + \theta(\lie{p}_1), z + v + \lie{p}_3' ) 
\qquad 
(v \in \lls),
\end{equation*}
which proves \eqref{eqn:intersection.of.tangents}.  

{From} \eqref{eqn:intersection.of.tangents}, we get
\begin{equation}
\label{eqn:intersection.subset.tangent.of.orbit}
T_{\xi} \calorbit \cap T_{\xi} \thetacalorbit_w \cap T_{\xi} X 
\subset \{ ( z + \lie{p}_1, z + \theta(\lie{p}_1), z + \lie{p}_3' ) \mid z \in \lie{k} \} = T_{\xi} \orbit .
\end{equation}
To see this, 
we concentrate on the third component of $ T_{\xi} X $, which must be of the form 
$ z' + \lie{p}_3' $ for some $ z' \in \lie{k} $.  
Equating this with the third component $ z + v + \lie{p}_3' $ of $ T_{\xi} \calorbit \cap T_{\xi} \thetacalorbit_w $, 
we get 
\begin{equation*}
(z - z') + v \in \lie{p}_3' 
\qquad
(z, z' \in \lie{k}; \;\; v \in \lie{s}) .
\end{equation*}
Since $ \lie{p}_3' $ is $ \theta $-stable, 
we get 
$ z - z' \in \llk \cap \llp_3' $ and $ v \in \lls \cap \llp_3' $.  
Therefore, the third component becomes $ z + v + \lie{p}_3' = z + \llp_3' $.

By Equations 
\eqref{eqn:tangent.of.orbit.subset.intersection} and 
\eqref{eqn:intersection.subset.tangent.of.orbit}, 
we have
\begin{equation*}
T_{\xi} \orbit \subset T_{\xi}( \calorbit \cap \thetacalorbit_w \cap X ) \subset T_{\xi} \calorbit \cap T_{\xi} \thetacalorbit_w \cap T_{\xi} X 
\subset T_{\xi} \orbit, 
\end{equation*}
and conclude that all the containments in the above formula are in fact equalities.  
This means that $ \orbit $ is an open neighborhood of $ \xi \in \calorbit \cap \thetacalorbit_w \cap X $.  
Since $ \xi $ is arbitrary, $ \calorbit \cap \thetacalorbit_w \cap X $ is smooth and 
its irreducible components coincide with connected components.  
Since the number of irreducible components of an algebraic variety is finite, 
we conclude that there are only finitely many $ K $-orbits in $ \calorbit \cap \thetacalorbit_w \cap X $.

Thus we finished the proof of Theorem \ref{theorem:finiteness.of.3.flags.implies.finiteness.of.2.flags}.
\end{proof}

The above theorem is strong enough to produce many interesting examples of double flag varieties of finite type.  
However, it also misses many possibilities.  
Here we introduce another kind of technique, which can present some more examples.  
A key idea is that a homogeneous space $G/Q$ sometimes can be embedded into
a product of (partial) flag varieties.
It is an equivariant compactification, and is considered to be a generalization of 
a complexification of the Harish-Chandra embedding
of a symmetric space into the product of the flag varieties.

First, let us explain the classical Harish-Chandra embedding.  

Let us assume that $ K $ is an intersection of a parabolic subgroup $ P $ of $ G $ and its opposite $ \opp{P} $.  
Thus $ K = P \cap \opp{P} $ is a Levi component of $ P $.  
Then $ G/K $ can be embedded into $ \GFl_{P} \times \GFl_{\opp{P}} $:
\begin{equation*}
G/K \ni g K \mapsto ( g P, g \opp{P} ) \in \GFl_{P} \times \GFl_{\opp{P}} ,
\end{equation*}
and this embedding is an open embedding (compare their dimension).  
Let us fix a Borel subgroup $ B \subset G $, and consider 
an embedding 
\begin{equation*}
B \backslash G/K  \inclusion B \backslash ( \GFl_{P} \times \GFl_{\opp{P}} ) \simeq 
G \backslash ( \GFl_B \times \GFl_{P} \times \GFl_{\opp{P}} )
\end{equation*}
Since $ \# B \backslash G / K < \infty $, there is an open $ B $-orbit, 
hence $ \GFl_{P} \times \GFl_{\opp{P}} $ is a spherical $ G $-variety.  
(See \S~\ref{subsection:spherical.varieties} for fundamental properties of spherical varieties.)  
Therefore there are finite number of $ B $-orbits on $ \GFl_{P} \times \GFl_{\opp{P}} $, 
which is equivalent to the finiteness of $ G $-orbits in 
$ \GFl_B \times \GFl_{P} \times \GFl_{\opp{P}} $.

Thus we proved the following

\begin{proposition}
\label{proposition:exist.open.orbit.via.HC.embedding}
Assume that a Levi component of a parabolic subgroup $ P $ is a symmetric subgroup of $ G $.  
Then 
$ \GFl_B \times \GFl_{P} \times \GFl_{\opp{P}} $ 
contains finitely many $ G $-orbits, where $ \opp{P} $ denotes a parabolic subgroup opposite to $ P $.
\end{proposition}

The assumption of the proposition above is satisfied for a symmetric pair $ (G, K) $ 
which is the complexification of a Hermitian symmetric pair $ (G_{\R}, K_{\R}) $.  
If $ P $ has an abelian unipotent radical, 
then $ K = P \cap \opp{P} $ satisfies this assumption (i.e., $ (G, K) $ is a symmetric pair; see \cite{Richardson.Roehrle.1992}).

Now we generalize the above situation to get a simple criterion of finiteness of $ K $-orbits on 
the double flag variety.

\begin{theorem}
\label{proposition:finiteness.by.HC.embedding}\label{theorem:finiteness.by.HC.embedding}
Let $ P^i \; (i = 1, 2, 3) $ be a parabolic subgroup of $ G $.  
Suppose that $ \GFl_{P^1} \times \GFl_{P^2} \times \GFl_{P^3} $
has finitely many $G$-orbits
and that $Q:= P^2 \cap P^3$ is a parabolic subgroup of $K$.
Then $\GKFl{P^1}{Q}$ has finitely many $K$-orbits.  

Moreover, if $ P^1 $ is a Borel subgroup $ B $ and the product $ P_2 P_3 $ is open in $ G $, 
then the converse is also true, i.e., 
the double flag variety 
$ \GKFl{B}{Q} $ 
is of finite type 
if and only if 
the triple flag variety 
$ \GFl_{B} \times \GFl_{P^2} \times \GFl_{P^3} $ 
is of finite type.
\end{theorem}

\begin{proof}
We have a $G$-equivariant diagonal embedding
$G/Q \inclusion \GFl_{P^2} \times \GFl_{P^3}$
by $ g \, Q \mapsto (g P^2, g P^3)$.  
Then we have the following natural inclusion
\begin{equation}
\label{eqn:HC.embedding.double.flag.variety}
K \backslash (\GKFl{P^1}{Q}) 
\cong P^1 \backslash G/Q
= G \backslash (\GFl_{P^1} \times G/Q)
\inclusion 
G \backslash (\GFl_{P^1} \times \GFl_{P^2} \times \GFl_{P^3}), 
\end{equation}
which proves the first claim.  

Let us assume that $ P^1 = B $ is a Borel subgroup and $ P_2 P_3 \subset G $ is open.  
To prove the converse, let us assume that $ \GKFl{B}{Q} $ is of finite type.  
Since $ K \backslash \GKFl{B}{Q} \simeq B \backslash G/Q $, 
there is an open $ B $-orbit on $ G/Q $.  
Since $ P_2 P_3 $ is open in $ G $, the map 
$G/Q \inclusion \GFl_{P^2} \times \GFl_{P^3}$ above is an open embedding, 
and consequently there is an open $ B $-orbit on $ \GFl_{P^2} \times \GFl_{P^3} $.  
Thus $ \GFl_{P^2} \times \GFl_{P^3} $ is a spherical $ G $-variety, hence 
there are only finitely many $ B $-orbits on it.  
Now, since 
$ G \backslash ( \GFl_{B} \times \GFl_{P^2} \times \GFl_{P^3} ) \simeq B \backslash ( \GFl_{P^2} \times \GFl_{P^3} ) $, 
we are done.
\end{proof}

Note that $ (P^1, P^2, P^3) = (B, P, \opp{P}) $ and $ Q = P \cap \opp{P} = K $ 
in Proposition \ref{proposition:exist.open.orbit.via.HC.embedding} above.

\section{Richardson-Springer theory}
\label{section:corollaries.of.main.theorem}

We use the same notation as in the former section.  
If $ P' = G $ and $ P = B $, Theorem \ref{theorem:finiteness.of.3.flags.implies.finiteness.of.2.flags} 
reduces to the one which claims that $ K \backslash G / B $ is a finite set.  
Let us compare this to the classification of $ K $-orbits on $ \GFl_B = G/B $ by Richardson-Springer.  

\subsection{Review of Richardson-Springer Theory}

First, we briefly review the theory of Richardson and Springer \cite{Richardson.Springer.1990, Richardson.Springer.1993}.  
We fix a $ \theta $-stable Borel subgroup $ B $ and a $ \theta $-stable maximal torus $ T \subset B $.  
Such pair always exists (\cite[Theorem~7.5]{Steinberg.1968}).
Let $ \maxT $ be the set of maximal tori in $ G $, 
and $ \maxTtheta $ the set of $ \theta $-stable maximal tori.  
As before $ \GFl_B $ denotes the set of all Borel subgroups in $ G $.  
We put $ \CTX = \{ ( T_1, B_1 ) \in \maxT \times \GFl_B \mid T_1 \subset B_1 \} $.  
Then there are natural projections $ p_1 : \CTX \to \maxT $ and $ p_2 : \CTX \to \GFl_B $.  

The projection $ p_2 : \CTX \to \GFl_B $ gives $ \CTX $ the structure of fiber bundle over $ \GFl_B $ with the fiber $ B_1/ T_1 $.  
The projection $ p_1 : \CTX \to \maxT $ is a Galois covering map with the Galois group $ W = W_G(T_1) $.  
Both of them tell us that $ \CTX $ is isomorphic to $ G/T $:
\begin{equation*}
\CTX \simeq G \times_B ( B / T ) \simeq G \times_{N_G(T)} ( N_G(T)/T ) \simeq G/T .
\end{equation*}
Put
\begin{equation*}
\CthetaTX := \{ ( T_1, B_1 ) \in \maxTtheta \times \GFl_B \mid T_1 \subset B_1 \} = \CTX \cap ( \maxTtheta \times \GFl_B ).
\end{equation*}

\begin{theorem}[Richardson-Springer]
\label{theorem:RS.KCtheta.isomorphic.to.KGB}
The $ K $-equivariant projection $ p_2 : \CthetaTX \to \GFl_B $ induces a bijection 
$ K \backslash \CthetaTX \xrightarrow{\;\sim\;} K \backslash \GFl_B $.
\end{theorem}

\begin{corollary}
Let us fix representatives $ \{ T_1 \} $ in the $ K $-conjugacy classes of the $ \theta $-stable maximal tori $ K \backslash \maxTtheta $.  
For each representative $ T_1 $, we also fix a Borel subgroup $ B_1 $ which contains $ T_1 $.  
Then there is a bijection
\begin{equation*}
\coprod_{T_1 \in K \backslash \maxTtheta} W_K(T_1) \backslash W_G(T_1)  \xrightarrow{\;\sim\;} K \backslash \GFl_B , 
\qquad
W_K(T_1) w \mapsto K \cdot (w B_1 w^{-1}) , 
\end{equation*}
where $ W_H(T_1) = N_H(T_1) / Z_H(T_1) $ is a Weyl group with representatives in $ H \subset G $.
\end{corollary}

The incidence variety $ \CthetaTX $ is sometimes too big for our purpose.  
We can take a smaller subvariety as follows.  
Define a map $ \tau : G \to G $ by $ \tau(g) = g^{-1} \theta(g) \;\; (g \in G) $.  
We denote by $ \Xi = \Image \tau $ the image of the map, which is known to be closed in $ G $ (\cite{Richardson.1982}).  
Since $ \tau $ is clearly invariant under the left translation by $ K $, it induces a map 
$ \Psi : K \backslash G \to \Xi $.  
By \cite[Lemma~2.4]{Richardson.1982}, $ \Psi $ is an isomorphism from the symmetric variety $ K \backslash G $ 
to the closed subvariety $ \Xi \subset G $.  

Recall the fixed $ \theta $-stable maximal torus $ T $.  
We define 
\begin{equation*}
\RSV := \tau^{-1}(N_G(T)) 
= \{ g \in G \mid \text{$ g^{-1} \theta(g) $ normalizes $ T $} \} , 
\end{equation*}
on which $ K $ acts on the left and $ T $ acts on the right.

\begin{theorem}[Richardson-Springer]
\label{theorem:RS.KVT.isomorphic.to.KGB}
There is a bijection $ K \backslash \RSV /T \xrightarrow{\;\sim\;} K \backslash \GFl_B $, which is induced by 
$ \RSV \ni g \mapsto g \cdot B \in \GFl_B $.
\end{theorem}

Let us briefly explain that this is an immediate consequence of Theorem \ref{theorem:RS.KCtheta.isomorphic.to.KGB}.
An element $ ( T_1, B_1 ) \in \CthetaTX $ is expressed as $ (T_1, B_1) = ( g T g^{-1}, g B g^{-1} ) $ for some $ g \in G $.  
The representative $ g \in G $ is determined up to the right multiplication of $ T $.  
Since $ T_1 = g T g^{-1} $ is $ \theta $-stable, we have
\begin{equation*}
g T g^{-1} = \theta(g T g^{-1}) = \theta(g) \, T \, \theta(g)^{-1} .
\end{equation*}
Hence $ g^{-1} \theta(g) \in N_G(T) $.  
Thus $ \RSV / T $ corresponds to $ \CthetaTX $ naturally by 
$ \RSV \ni g \mapsto ( g \cdot T, g \cdot B ) \in \CthetaTX $.  
So $ K $-orbits in $ \RSV / T $ are in bijection with $ K $-orbits in $ \CthetaTX $, 
which are in turn bijective to $ K \backslash \GFl_B $.

Now we get a map 
\begin{equation*}
K \bsl \GFl_B \xrightarrow{\;\sim\;} K \bsl \RSV / T \xrightarrow{\;\tau\;} N_G(T) \to W = N_G(T)/T , 
\end{equation*}
which sends $ K g B $ to $ w = g^{-1} \theta(g) $ in $ W $.  
Note that $ \theta(w) = \theta(g)^{-1} g = w^{-1} $.  
We call $ v \in W $ a twisted involution if $ \theta(v) = v^{-1} $ holds, and put
$ \TwistedInv = \{ v \in W \mid \theta(v) = v^{-1} \} $, the set of twisted involutions.  
With this notation, we finally get a map
\begin{equation}
\label{eqn:map.from.KGB.to.I}
\phi : K \bsl \GFl_B \to \TwistedInv \subset W , 
\end{equation}
which we call the \emph{Richardson-Springer map}.

\subsection{Geometry of Richardson-Springer map}

Let us recall the notation in \S~\ref{section:double.flags.for.s.pair}.  
We take $ P = B $, a Borel subgroup, and $ P' = G $ so that 
$ \GKFl{P}{Q} = \GFl_B $.  

We take a $ G $-orbit $ \calorbit $ in $ \GFl_B \times \GFl_B $ under the diagonal action.  
Since $ G \bsl ( \GFl_B \times \GFl_B ) \simeq B \bsl G /B $, the $ G $-orbits are classified by the Weyl group $ W = W_G(T) $.  
So we write $ \calorbit = \calorbit_w \; (w \in W) $.  
Let us consider the $ \theta $-twisted embedding of $ \GFl_B $ into $ \GFl_B \times \GFl_B $, i.e.,
\begin{equation*}
\Dtheta : \GFl_B \injection \GFl_B \times \GFl_B , 
\quad
B_1 \mapsto ( B_1, \theta(B_1) ) .
\end{equation*}
We denote $ X = \Dtheta(\GFl_B) $.  
Then Lemma \ref{lemma:finite.orbits.on.O.cap.Otheta} tells us the following

\begin{lemma}
For each $ w \in W $, 
the connected components of $ \calorbit_w \cap X $ are 
exactly the irreducible components.  
Each connected component is a $ K $-orbit, hence 
there are finitely many $ K $-orbits in $ \calorbit_w \cap X $.
\end{lemma}

Now pick a point $ \xi $ in $ \calorbit_w \cap X $.  
Then $ \xi = ( B_1 , \theta(B_1) ) = (g \cdot B, (g \dot{w}) \cdot B) $, 
where $ \dot{w} \in N_G(T) $ represents $ w \in W = N_G(T)/T $.  
\begin{equation*}
B_1 = g B g^{-1} , \quad
\theta(B_1) = (g \dot{w})  B (g \dot{w})^{-1}, \quad
(g \dot{w})^{-1} \theta(g) \in B
\end{equation*}
Thus we have $ \dot{w}^{-1} g^{-1} \theta(g) \in B $.  From Theorem \ref{theorem:RS.KVT.isomorphic.to.KGB}, 
we can assume $ g^{-1} \theta(g) \in N_G(T) $.  
Therefore $ \dot{w}^{-1} g^{-1} \theta(g) \in B \cap N_G(T) = T $.  
Thus $ g^{-1} \theta(g) $ represents $ w $ also.

\begin{theorem}
Let us denote $ X = \Dtheta(\GFl_B) \subset \GFl_B \times \GFl_B $.  
For $ w \in W $, let us consider a $ G $-orbit $ \calorbit_w = G \cdot (B, w \cdot B) \in \GFl_B \times \GFl_B $.  
If $ \calorbit_w \cap X \neq \emptyset $, then $ w \in \TwistedInv $ is a twisted involution, i.e., it satisfies $ w^{-1} = \theta(w) $.  
Moreover, if $ w \in \TwistedInv $, 
the connected components of $ \calorbit_w \cap X $ correspond bijectively to the $ K $-orbits in $ K \bsl \GFl_B $ 
which are in the fiber $ \phi^{-1}(w) $ 
of $ w $ of the map $ \phi $ 
{\upshape(}see Equation \eqref{eqn:map.from.KGB.to.I}{\upshape)}.
\end{theorem}

This theorem gives a geometric interpretation of the Richardson-Springer map $ \phi : K \bsl \GFl_B \to \TwistedInv $.

\section{Spherical actions on multiple flag varieties}
\label{section:spherical.action}

\subsection{Spherical varieties}
\label{subsection:spherical.varieties}

The finiteness of $ K $-orbits on the product of flag varieties and 
spherical actions of $ G $ or $ K $ are closely related.  

Recall that a $ G $-variety $ X $ is called spherical if it has an open dense $ B $-orbit, 
where $ B $ is a Borel subgroup.  
Note that $ X $ is $ G $-spherical if and only if $ \# B \backslash X < \infty $.

Let us begin with an easy but fundamental lemma.

\begin{lemma}
\label{lemma:sphericality.and.double.flags.of.finite.type}
The following conditions are equivalent.
\begin{thmenumerate}
\item 
Let $B \times S$ be a Borel subgroup of $G \times K$.
Then $K$ has finitely many orbits on $\GKFl{B}{S}$.
\item 
$G \cong (G \times K)/ \Delta K$ is $G\times K$-spherical.
\item Every irreducible finite-dimensional rational representation of $G$
is decomposed into the representations of $K$ multiplicity freely.
\end{thmenumerate}
\end{lemma}

\begin{proof}
Let us consider the condition (1).  
Since $ \GKFl{B}{S} = (G \times K) / (B \times S) $, 
the finiteness of $ K $-orbits on it implies 
the finiteness of $ B \times S $-orbits on $ K \backslash (G \times K) $.  
Since $ B \times S $ is Borel in $ G \times K $, 
this is equivalent to that $ K \backslash (G \times K) $ is 
$ (G \times K) $-spherical, which is the condition (2).  

Note that $ (G \times K)/K $ is affine.  
So, the existence of an open $ B \times S $-orbit is equivalent to the condition that 
the regular function ring $ \C[G \times K]^K $ decomposes multiplicity freely 
as a representation of $ G \times K $.  
By the Frobenius reciprocity, we know
\begin{equation*}
\C[G \times K]^K \simeq \bigoplus_{(\pi, \tau) \in \irreps{G} \times \irreps{K}} 
\Hom_K(\pi, \tau) \otimes (\pi \boxtimes \tau^{\ast}) ,
\end{equation*}
where $ \tau^{\ast} $ denotes the contragredient of $ \tau $ and 
$ \pi \boxtimes \tau^{\ast} $ means outer tensor product.  
Thus we get $ \dim \Hom_K(\pi, \tau) \leq 1 $ for any $ \pi \in \irreps{G} $ and $ \tau \in \irreps{K} $, 
which is equivalent to the condition (3).
\end{proof}

There are very few examples which satisfy the conditions in the above lemma, 
and they are completely classified by Kr\"{a}mer \cite{Kramer.1976}.    
Best known one might be the pair $ (G, K) = (\GL_n, \GL_1 \times \GL_{n - 1}) $, which 
is related to the multiplicity free branching rule (Pieri formula) and Gelfand-Zeitlin basis.

The following theorem is a direct consequence of 
Theorem \ref{theorem:finiteness.of.3.flags.implies.finiteness.of.2.flags}.

\begin{theorem}
\label{cor:spherical.2.flags.implies.K.spherical.flags}\label{theorem:Kspherical.partial.flags}
Let $ P $ be a parabolic subgroup of $ G $.  
If $ \GFl_P \times \GFltheta_P $ is a spherical $ G $-variety, then 
$ \GFl_P $ is a spherical $ K $-variety.
\end{theorem}

\begin{proof}
Since $ G \backslash ( \GFl_P \times \GFltheta_P \times \GFl_B ) \simeq B \backslash ( \GFl_P \times \GFltheta_P ) $, 
the product $ \GFl_P \times \GFltheta_P $ is a spherical $ G $-variety if and only if 
there are finitely many $ G $-orbits on $ \GFl_P \times \GFltheta_P \times \GFl_B $.  
We can assume that $ B $ is $ \theta $-stable and $ S := K \cap B $ is a Borel subgroup of $ K $.  
Then, by Theorem \ref{theorem:finiteness.of.3.flags.implies.finiteness.of.2.flags}, 
this implies that there are finitely many $ K $-orbits on $ \GKFl{P}{S} \simeq G/P \times K/S $.  
Since $ K \backslash ( G/P \times K/S ) \simeq S \backslash G / P $, 
this is equivalent to say that $ \GFl_P \simeq G/P $ is $ K $-spherical.
\end{proof}

Let us take a $ \theta $-stable Borel subgroup $ B $ of $ G $, and fix a positive root system $ \Delta^+ $ 
corresponding to $ B $.  We denote by $ \Pi \subset \Delta^+ $ a simple system.  
For a parabolic subgroup $ P $ in $ G $ which contains $ B $, we can 
associate a subset $ \Phi \subset \Pi $ so that 
$ \Pi \setminus \Phi $ generates a sub root system for a Levi component of $ P $.  
For $ \alpha \in \Phi $, we denote a fundamental weight corresponding to $ \alpha $ by $ \omega_{\alpha} $.  
Put $ \lambda = \sum_{\alpha \in \Phi} c_{\alpha} \omega_{\alpha} $ 
a linear combination of those fundamental weights with positive integer coefficients $ c_{\alpha} $'s.
We assume that $ \lambda $ is integral for $ G $.  
Let us denote by $ V_{\lambda} $ a finite dimensional irreducible representation of $ G $ 
with highest weight $ \lambda $, and by $ v_{\lambda} $ its highest weight vector.  
Then we have 
\begin{equation*}
P = \{ g \in G \mid g \cdot v_{\lambda} \in \C v_{\lambda} \}  .
\end{equation*}
If we denote by $ \bbP(V_{\lambda}) $ a projective space over $ V_{\lambda} $ and 
$ [v_{\lambda}] $ a point in $ \bbP(V_{\lambda}) $ determined by the line through it, 
it is equivalent to say that 
$ G \cdot [v_{\lambda}] \simeq \GFl_P $.  
Let us denote $ \affGFl_P = \closure{G v_{\lambda}} \subset V_{\lambda} $, 
the affine cone over $ \GFl_P $.  

With these notations, we have the following

\begin{lemma}
\label{lemma:K.spherical.PFV.and.K.MF}
The partial flag variety $ \GFl_P $ is  $ K $-spherical 
if and only if 
$ V_{k \lambda}^{\ast} \restrict_K $ decomposes multiplicity freely as a $ K $-module 
for any non-negative integer $ k \geq 0 $.
\end{lemma}

\begin{proof}
The partial flag variety $ \GFl_P $ is  $ K $-spherical 
if and only if the affine cone $ \affGFl_P $ is $ \C^{\times} \times K $-spherical.  
Since $ \affGFl_P $ is an affine variety, 
it is $ \C^{\times} \times K $-spherical 
if and only if 
the regular function ring 
$ \C[\affGFl_P] $ decomposes multiplicity freely.
Note that 
\begin{equation*}
\C[\affGFl_P] \simeq \bigoplus\nolimits_{k \geq 0} V_{k \lambda}^{\ast}
\end{equation*}
as a $ G $-module.  
Since $ \C^{\times} $-action specifies one of $ V_{k \lambda}^{\ast} $, 
the restriction $ V_{k \lambda}^{\ast} \restrict_K $ is a multiplicity free $ K $-module.
\end{proof}

Without loss of generality, we can assume that the root system $ \Delta $ is defined with respect to 
a $ \theta $-stable maximal torus $ T $.  
Thus there is a well-defined action of $ \theta $ on the root system $ \Delta $.  
Since $ B $ is assumed also to be $ \theta $-stable, 
$ \theta $ preserves the simple system $ \Pi $, 
and we easily see that $ \theta(P) $ corresponds to $ \theta(\Phi) $.  
Put $ \lambda^{\theta} = \sum_{\alpha \in \theta(\Phi)} c_{\alpha} \omega_{\alpha} $.  
Under an obvious notation, we conclude that 

\begin{lemma}
The product of partial flag varieties $ \GFl_P \times \GFltheta_P $ is  $ G $-spherical 
if and only if 
$ V_{k \lambda} \otimes V_{\ell \lambda^{\theta}} $ decomposes multiplicity freely as a $ G $-module 
for any non-negative integers $ k , \, \ell \geq 0 $.
\end{lemma}

The proof is the same as Lemma \ref{lemma:K.spherical.PFV.and.K.MF} 
(and essentially, this follows from the lemma if we consider $ \bbG = G \times G $ and $ \bbK = \Delta (G) $).

Thus we can reinterpret Theorem \ref{theorem:Kspherical.partial.flags} as follows.

\begin{corollary}
Let $ P $ be a parabolic subgroup containing a $ \theta $-stable Borel subgroup $ B $ of $ G $, 
and we assume the notations above.  
If the tensor product $ V_{k \lambda} \otimes V_{\ell \lambda^{\theta}} $ is a multiplicity free $ G $-module 
for any non-negative integers $ k , \, \ell \geq 0 $, 
then the restriction 
$ V_{m \lambda} \restrict_K $ decomposes multiplicity freely as a $ K $-module for any $ m \geq 0 $.
\end{corollary}

\subsection{Maximally split parabolic in a real form}

Let $ \lier{g} $ be a real Lie algebra which is a real form of $ \lie{g} $.  
Let $ G_{\R} $ be a connected analytic Lie subgroup in $ G $ corresponding to $ \lier{g} $, and 
we assume it is non-compact.  
Choose a maximal compact subgroup $ K_{\R} $ of $ G_{\R} $.  
Then we have a Cartan decomposition 
$ \lier{g} = \lier{k} \oplus \lier{s} $ 
corresponding to $ K_{\R} $.  
It is well known that, 
for a symmetric pair $ (G, K) $, there always exists 
such a non-compact Riemannian symmetric pair $ (G_{\R}, K_{\R}) $, and 
our involution $ \theta $ coincides with the complexification of the Cartan involution 
associated to $ G_{\R}/K_{\R} $.  

Choose a maximal abelian subspace $ \lier{a} $ in $ \lier{s} $.  
Then a choice of a positive system of the restricted root system 
$ \Sigma(\lier{g}, \lier{a}) $ determines a real parabolic subalgebra $ \lier{p} $ which is maximally split in $ \lier{g} $.  
Let $ \lie{p}_{\mathrm{min}} $ be the complexification of $ \lier{p} $, and $ \miniP $ 
the corresponding complex parabolic subgroup of $ G $.  
We denote by $ \GFl_{\miniP} \simeq G/\miniP $ a partial flag variety of parabolic subgroups conjugate to $ \miniP $ as usual.  

\begin{lemma}
\label{lemma:dense.Korbit.simeq.KoverM}
The dense open $ K $-orbit in $ \GFl_{\miniP} $ is isomorphic to $ K/M $, 
where $ M = Z_{K}(\lie{a}) $ is the centralizer of $ \lie{a} = \C \otimes_{\R} \lier{a} $ in $ K $.
\end{lemma}

\begin{proof}
This lemma is well known, but we prove it for the sake of self-containedness.  
Since $ K \cap \miniP = M $, the $ K $-orbit through $ \lie{p}_{\mathrm{min}} \in \GFl_P $ 
is isomorphic to $ K/M $.  
The complex dimension of $ K/M $ is equal to the real dimension of Iwasawa's $ \lier{n} $, 
which is also equal to the dimension of $ \GFl_{\miniP} $.  
Thus, the orbit must be an open orbit.
\end{proof}

The following is a corollary to Theorem \ref{cor:spherical.2.flags.implies.K.spherical.flags}.

\begin{corollary}
Let $ \miniP $ be the complexification of a maximally split parabolic subgroup of $ G_{\R} $ as above.  
If $ \GFl_{\miniP} \times \GFl_{\miniP} $ is a $ G $-spherical variety, 
then $ K/M $ is a spherical $ K $-variety.
\end{corollary}

\begin{proof}
If $ \GFl_{\miniP} \times \GFl_{\miniP} $ is a $ G $-spherical variety, 
then $ \GFl_{\miniP} $ is a spherical $ K $-variety by 
Theorem \ref{cor:spherical.2.flags.implies.K.spherical.flags} 
(note that we can choose a $ \theta $-stable parabolic from $ \GFl_{\miniP} $ 
so that $ \GFl_{\miniP} = \GFltheta_{\miniP} $).  
So there are finitely many $ S $-orbits, where $ S $ is a Borel subgroup of $ K $.  
Since $ K/M $ is an open orbit in $ \GFl_{\miniP} $, there are only finitely many $ S $-orbits in $ K/M $, 
which implies $ K/M $ is $ K $-spherical.
\end{proof}

We have a partial converse to the above corollary.

\begin{proposition}
The $ K $-variety $ K/M $ is spherical 
if and only if $ \GKFl{\miniP}{Q} $ contains finitely many $ K $-orbits for any parabolic subgroup $ Q $ of $ K $.
\end{proposition}

Note that if $ (K, M) $ is a symmetric pair, then $ K/M $ is spherical.

\begin{proof}
Let $ B \subset G $ be a $ \theta $-stable Borel subgroup such that $ S := K \cap B $ is a Borel subgroup of $ K $.  
We can consider the Borel subgroup $ S $ instead of $ Q $ without loss of generality.  
The variety $ K/M $ is $ K $-spherical if and only if there exists an open $ S $-orbit on $ K/M $.  
By Lemma \ref{lemma:dense.Korbit.simeq.KoverM}, $ K/M $ is an open $ K $-orbit in $ \GFl_{\miniP} $.  
So the open $ S $-orbit in $ K/M $ turns out to be an open $ S $-orbit in $ \GFl_{\miniP} $, 
which means $ \GFl_{\miniP} $ is a spherical $ K $-variety.  
Thus we have finitely many $ S $-orbits on $ \GFl_{\miniP} $.  
Through the isomorphism $ K \backslash ( \GKFl{\miniP}{S} ) \simeq S \backslash \GFl_{\miniP} $, 
we conclude that $ \GKFl{\miniP}{S} $ also contains finitely many $ K $-orbits.
\end{proof}

\section{Double flag varieties of type A}
\label{section:typeA}

Let us consider a group $ G $ of type A.  
There are three types of symmetric pairs $ (G, K) $, denoted by AI, AII, AIII (see \cite[\S~X.6]{Helgason.DG.1978}).  
Namely they are 
$ \SL_n / \SO_n $, $ \SL_{2m} / \Sp_{2m} $, and 
$ \GL_n/ \GL_p \times \GL_q \;\; (n = p + q) $.  
We will construct examples of double flag varieties with finitely many $ K $-orbits, using 
Theorem \ref{theorem:finiteness.of.3.flags.implies.finiteness.of.2.flags}.  

Recall the notation $ P_{\lambda} $ for an (upper triangular) standard parabolic subgroup of $ \GL_n $ from \S \ref{subsection:MWZ.typeA}, 
where $ \lambda = ( \lambda_1, \dots, \lambda_{\ell}) $ is a composition of size $ n $.  
In fact, $ P_{\lambda} $ is realized as the stabilizer of a partial flag of subspaces in $ \C^n $ 
of dimension $ \lambda_1, \lambda_1 + \lambda_2, \dots, \lambda_1 + \cdots + \lambda_{\ell} $.

\subsection{Type AI and AII}

Let $ G/K = \SL_n / \SO_n \; ( n \geq 3 ) $ or 
$ G/K = \SL_{2m} / \Sp_{2m} \; (m \geq 2) $.  
In these cases, a mirabolic parabolic subgroup is not conjugate to a $ \theta $-stable parabolic subgroup.  
So we have less possibilities to apply 
Theorem \ref{theorem:finiteness.of.3.flags.implies.finiteness.of.2.flags}.

\begin{proposition}
Let $ (G, K) = (\SL_n, \SO_n) $ or 
$ (G, K) = (\SL_{2m}, \Sp_{2m}) $, which is  a symmetric pair of type {\upshape AI} or {\upshape AII} respectively. 
If $ P \subset G $ and $ Q \subset K $ are a pair of parabolic subgroups among the following list {\upshape(1)--(2)}, 
then there are finitely many $ K $-orbits in $ \GKFl{P}{Q} $.  
\begin{thmenumerate}
\item
$ P $ is a maximal parabolic subgroup of $ G $, and $ Q $ is an arbitrary parabolic subgroup of $ K $.
\item
Assume that $ n \geq 4 $ is an even integer if 
$ (G, K) = (\SL_n, \SO_n) $.  
$ P = P_{\lambda} $ is a parabolic subgroup of $ G $ 
with $ \lengthof{\lambda} = 3 $ {\upshape (}i.e., $ \lambda = ( \lambda_1, \lambda_2 , \lambda_3 ) ${\upshape )}, 
and $ Q $ is a Siegel parabolic subgroup of $ K $.  
Here we say $ Q $ is a Siegel parabolic subgroup if it is the stabilizer of a maximal isotropic space.  
\end{thmenumerate}
\end{proposition}

\begin{proof}
Here we only give a proof for type AI.  
The proof for type AII is similar.

(1)\  
Type $ D_{r + 2} $ in Theorem \ref{theorem:MWZtypeA} implies the result.

(2)\ 
Put $ n = 2m $.  
We use type $ E_6 $ in Theorem \ref{theorem:MWZtypeA}.  
Since the maximal parabolic $ P' $ in the list should be $ \theta $-stable, 
it must be a parabolic subgroup of $ \SL_{2m} $ corresponding to a partition $ (m, m) $.  
So we can take $ Q = P' \cap K $ as a Siegel parabolic of $ \SO_{2m} $ with an appropriate choice of conjugates of $ P' $.
\end{proof}

\subsection{Type AIII}
\label{subsection:typeAIII.general.case}

$ G/K = \GL_n/ \GL_p \times \GL_q \;\; (n = p + q) $.  


We get three cases in which the double flag variety $ \GKFl{P}{Q} $ has finitely many $ K $-orbits.  
These are direct consequence of 
Theorem \ref{theorem:finiteness.of.3.flags.implies.finiteness.of.2.flags}.  

\begin{proposition}
Let $ G/K = \GL_n/ \GL_p \times \GL_q $ be a symmetric space of type {\upshape{}AIII}.
Let $ P $ be a parabolic subgroup of $ G $ and $ Q $ that of $ K $.  
If $ P $ and $ Q $ are among the following list {\upshape(1)--(3)}, 
then there are finitely many $ K $-orbits in $ \GKFl{P}{Q} $.
\begin{thmenumerate}
\item
$ P $ is any parabolic subgroup of $ G $, and 
$ Q = Q_1 \times Q_2 $ is a parabolic subgroup of $ K $ which satisfies
{\upshape(i)}\ 
$ Q_1 $ is of partition type $ (1, p -1) $ and $ Q_2 = \GL_q $; or
{\upshape(ii)}\ 
$ Q_1 = \GL_p $ and $ Q_2 $ is of partition type $ (q - 1, 1) $. 
\item
$ P $ is a maximal parabolic subgroup of $ G $, and 
$ Q $ is any parabolic in $ K $.
\item
$ P = P_{\lambda} $ is a parabolic subgroup of $ G $ which corresponds to a composition $ \lambda $ 
with $ \lengthof{\lambda} = 3 $ 
and $ Q $ is a maximal parabolic subgroup of $ K $.
\end{thmenumerate}
\end{proposition}

\begin{proof}
For (1), 
We use type $ S_{q, r} $ in Theorem \ref{theorem:MWZtypeA}.  
For (2), 
we use type $ D_{r + 2} $ in the same theorem, and 
for (3), we use type $ E_6 $.
\end{proof}

Few remarks are in order.

In Case (1) in the above theorem, 
$ \KFl_{Q} $ is isomorphic to a projective space $ \bbP(\C^p) $ or $ \bbP(\C^q) $.  
We call these double flag varieties ``mirabolic'' after 
\cite{Travkin.2009} and \cite{FGT.2009}.

In Case (2), if $ P = P_{(m, n -m)} $, 
$ \GFl_{P} $ is a Grassmannian $ \Grass_m(\C^n) $ of $ m $-dimensional subspaces in $ \C^n $ 
($ n = p + q $).  
Thus the action of $ K = \GL_p \times \GL_q $ on 
$ \Grass_m(\C^{p + q}) \times \GFl_{Q^1}^{\GL_p} \times \GFl_{Q^2}^{\GL_q} $ has 
finitely many orbits with obvious notations.

\medskip

If $ P = B $ is a Borel subgroup of $ G $, 
we are able to give a complete classification of 
the double flag variety $ \GKFl{B}{Q} $ of finite type 
for a symmetric pair of type AIII.  

\begin{theorem}
Let $G=\GL_n$ and 
$B \subset G$ be a Borel subgroup; 
$K=\GL_p \times \GL_q$ with $p+q=n$, 
$q \ge p \ge 1$;
and $Q^1$ is a parabolic subgroup of $\GL_p$ 
and $Q^2$ is that of  $\GL_q$.
Put $Q=Q^1  \times Q^2$ a parabolic subgroup of $K$.  
Then there are only finitely many $K$-orbits on $\GFl_B \times \KFl_Q$
if and only if 
$Q^1$ and $Q^2$ are in the following table.
{
\newcommand{\any}{\text{\upshape arbitrary}}
\newcommand{\maximal}{\text{\upshape maximal}}
\label{lemma.conditionC3:item.list.fo.Q2}
\[
\begin{array}{|c||c|c|c|c|}
\hline
\mbox{\rm Case} & p & Q^1 & Q^2 & \KFl_Q \\
\hline \hline
\mbox{\rm (i)} & \any & \GL_p & \GL_q & \{ \mbox{\upshape{point}} \} \\
\hline
\mbox{\rm (ii)} & \any & \GL_p & \mbox{\rm mirabolic} & \bbP(\C^q) \\
\hline
\mbox{\rm (iii)} & 1 & \GL_1 & \any  & \GL_q/Q^2 \\
\hline
\mbox{\rm (iv)} & 2 & \GL_2 & \maximal & \Grass_m(\C^q) \\
\hline
\mbox{\rm (v)} & \any &  \mbox{\rm mirabolic} & \GL_q & \bbP(\C^p) \\
\hline
\end{array}
\]
Here the second column indicates the condition on $p$.
}
\end{theorem}

\begin{proof}
We use Theorem \ref{proposition:finiteness.by.HC.embedding}.  
Let $\lambda$ be a composition of $p$,
and $\mu$ be that of $q$.
Note that $(\lambda,\mu)$ is a composition of $n$.
We put $P^2 = P_{(\lambda,\mu)}$, which is a standard parabolic subgroup of $ G $, and 
$P^3 = \opp{P_{(p,q)}} $, a parabolic subgroup of $G$ opposite to the standard parabolic $ P_{(p,q)} $.
It is easy to check that 
$Q=P^2 \cap P^3$ is a parabolic subgroup $P_\lambda \times P_\mu $ of $K=\GL_p \times GL_q$ 
and the product $ P^2 P^3 $ is open dense in $ G $.  
Note that $ P_{\lambda} $ (respectively $ P_{\mu} $) is a parabolic subgroup of $ \GL_p $ (respectively $ \GL_q $).  
Now we are in the setting of Theorem \ref{proposition:finiteness.by.HC.embedding}, and conclude that 
$ \GKFl{B}{Q} $ is of finite type if and only if 
the triple flag $ \GFl_{B} \times \GFl_{P_{(\lambda, \mu)}} \times \GFl_{P_{(p,q)}} $ is of finite type.  
From Theorem \ref{theorem:MWZtypeA}, we deduce the table above.
\end{proof}

\subsection{Summary}

As a summary, we give tables of the double flag varieties with finitely many $ K $-orbits in 
Tables \ref{table:typeAI}--\ref{table:AIII.general} below.  
Note that these tables do not exhaust all the cases.

\begin{table}[htbp]
\caption{Type AI : $ G/K = \SL_n / \SO_n \;\; ( n \geq 3 ) $}
\label{table:typeAI}
\medskip
{
\newcommand{\any}{\text{arbitrary}}
\newcommand{\maximal}{\text{maximal}}
\hfil
\begin{math}
\begin{array}{c|c|c|c|c}
\hline
P & Q & \GFl_{P} & \KFl_{Q} & \text{extra condition} \\
\hline
\maximal & \any & \Grass_m(\C^n) & \KFl_{Q} \\
(\lambda_1, \lambda_2, \lambda_3) & \text{Siegel} & \GFl_{P} & \LGrass(\C^n) & \text{$ n $ is even} \\
\hline
\end{array}
\end{math}
\hfil
%
\\[3ex]
%
\caption{Type AII : $ G/K = \SL_{2n} / \Sp_{2n} \;\; (n \geq 2) $}
\label{table:typeAII}
\medskip
\hfil
\begin{math}
\begin{array}{c|c|c|c}
\hline
P & Q & \GFl_{P} & \KFl_{Q} \\
\hline
\maximal & \any & \Grass_m(\C^{2n}) & \KFl_{Q} \\
(\lambda_1, \lambda_2, \lambda_3) & \text{Siegel} & \GFl_{P} & \LGrass(\C^{2n}) \\
\hline
\end{array}
\end{math}
\hfil
%
\\[3ex]
%
\caption{Type AIII : $ G/K = \GL_n/ \GL_p \times \GL_q \;\; (n = p + q) $. } 
\label{table:AIII.general}
\medskip
\hfil
\newcommand{\mirabolic}{\text{mirabolic}}
\begin{math}
\begin{array}{c|cc|c|c}
\hline
P & Q_1 & Q_2 & \GFl_{P} & \KFl_{Q} \\
\hline
\any & \mirabolic & \GL_q & \GFl_{P} & \bbP(\C^p) \\
\any & \GL_p & \mirabolic & \GFl_{P} & \bbP(\C^q) \\
\maximal & {\any} & {\any} & \Grass_m(\C^n) & \KFl_{Q} \\
(\lambda_1, \lambda_2, \lambda_3) & \GL_p & {\maximal} & \GFl_{P} & \Grass_k(\C^q) \\
(\lambda_1, \lambda_2, \lambda_3) & {\maximal} & \GL_q & \GFl_{P} & \Grass_k(\C^p) \\
\hline
\rule{0pt}{2.5ex}
{\any} & \GL_1 \; (p = 1) & {\any} & \GFl_{P} & \GL_q/Q_2 \\
{\any} & \GL_2 \; (p = 2) & {\maximal} & \GFl_{P} & \Grass_m(\C^q) \\
\hline
\end{array}
\end{math}
\hfil
}
\end{table}

\section{Double flag varieties of type C}
\label{section:typeC}

Let us consider a symmetric pair of type C.  
There are two irreducible symmetric spaces of type C; namely, type CI and CII.
So we consider a symmetric space $ G/K = \Sp_{2n}/ \GL_n $ of type CI, or 
$ G/K = \Sp_{2n} / \Sp_{2p} \times \Sp_{2q} $ of type CII ($n = p + q$) in this section.  

First, recall the notation $ P_{\lambda} $ of a standard parabolic subgroup of $ \Sp_{2n} $ from \S \ref{subsection:MWZ.typeC}, 
where $ \lambda $ is a composition of size $ 2 n $.  
The parabolic $ P_{\lambda} $ is realized as the stabilizer of a partial flag of isotropic subspaces in $ \C^{2n} $.  
In particular, a maximal parabolic subgroup $ P_{(m, 2 n - 2 m, m)} \; (0 < m \leq n) $ is the stabilizer of an isotropic subspace 
of dimension $ m $.  

If $ m = n $, then a totally isotropic subspace of dimension $ n $ is called \emph{Lagrangian}, and 
we denote by $ \LGrass(\C^{2n}) $ the Grassmannian of all the Lagrangian subspaces in $ \C^{2n} $.    
Let $ P_{(n, n)} $ be a Siegel parabolic subgroup, which fixes a Lagrangian subspace.  
Since $ G = \Sp_{2n} $ acts on $ \LGrass(\C^{2n}) $ transitively, 
we have $ G/P_{(n,n)} \simeq \LGrass(\C^{2n}) $.  

If $ m < n $, let $ \IGrass_m(\C^{2n}) $ be the Grassmannian of isotropic subspaces of fixed dimension $ m $.  
As in the case of the Lagrangian Grassmannian, 
we can identify $ G/P_{(m, 2n - 2m, m)} \simeq \IGrass_m(\C^{2n}) $.  
Note that, if $ m = 1 $, this reduces to $ G/P_{(1, 2n - 2, 1)} \simeq \bbP(\C^{2n}) $.

Theorem \ref{theorem:finiteness.of.3.flags.implies.finiteness.of.2.flags} gives us 
several examples of double flag varieties of finite type.

\begin{proposition}
Let $ G/K = \Sp_{2n}/ \GL_n $ be a symmetric space of type {\upshape{CI}} or 
$ G/K = \Sp_{2n} / \Sp_{2p} \times \Sp_{2q} $ of type {\upshape{CII}} $ (n = p + q) $.  
If a pair of parabolic subgroups $ P \subset G $ and $ Q \subset K $ is among the following list 
{\upshape(1)--(3)}, 
then the double flag variety $ \GKFl{P}{Q} $ is of finite type.
\begin{thmenumerate}
\item
$ P = P_{(n, n)} $ is a Siegel parabolic subgroup of $ G $, and $ Q \subset K $ is an arbitrary parabolic subgroup.
\item
$ P = P_{(1, 2 n - 2, 1)} $ is a maximal parabolic subgroup of $ G $, and $ Q \subset K $ is an arbitrary parabolic subgroup.
\item
Let us assume that $ G/K = \Sp_{2n} / \Sp_{2p} \times \Sp_{2q} $ is of type {\upshape{CII}}.  
$ P = P_{(m, 2 n - 2 m, m)} $ $ (1 < m < n) $ is a maximal parabolic subgroup of $ G $, 
and $ Q \subset K $ is a product of Siegel parabolic subgroups in $ \Sp_{2p} $ and $ \Sp_{2q} $.
\end{thmenumerate}
\end{proposition}

\begin{proof}
For (1), we use type $ \Sp{D_{r + 2}} $ in Theorem \ref{theorem:MWZtypeC} and apply 
Theorem \ref{theorem:finiteness.of.3.flags.implies.finiteness.of.2.flags}.  
Similarly, for (2), we use type $ \Sp{Y_{4, r}} $ in Theorem \ref{theorem:MWZtypeC}, and 
for (3), we use $ \Sp{E_6} $.
\end{proof}

As a summary, we give tables of the double flag varieties of type C with finitely many $ K $-orbits in 
Tables \ref{table:typeCI}--\ref{table:typeCII} below.  
Note that these tables do not exhaust all the cases.

\begin{table}[htbp]
\caption{Type CI : $ G/K = \Sp_{2n} / \GL_n \;\; ( n \geq 2 ) $}
\label{table:typeCI}
\medskip
{
\hfil
\newcommand{\any}{\text{any}}
\newcommand{\maximal}{\text{maximal}}
\newcommand{\Siegel}{\text{Siegel}}
\begin{math}
\begin{array}{c|c|c|c}
\hline
P & Q & \GFl_{P} & \KFl_{Q} \\
\hline
\Siegel & \any & \LGrass(\C^{2n}) & \KFl_{Q} \\
(1, 2 n - 2, 1) & \any & \bbP(\C^{2n}) & \KFl_{Q} \\
\hline
\end{array}
\end{math}
\hfil
}
%
\\[3ex]
%
\caption{Type CII : $ G/K = \Sp_{2n} / \Sp_{2p} \times \Sp_{2q} \;\; (n = p + q) $}
\label{table:typeCII}
\medskip
{
\hfil
\newcommand{\any}{\text{any}}
\newcommand{\maximal}{\text{maximal}}
\newcommand{\Siegel}{\text{Siegel}}
\begin{math}
\begin{array}{c|c|c|c}
\hline
P & Q & \GFl_{P} & \KFl_{Q}\\
\hline
\Siegel & \any & \LGrass(\C^{2n}) & \KFl_{Q} \\
(m, 2 n - 2m, m) & \Siegel \times \Siegel & \IGrass_m(\C^{2n}) & \LGrass(\C^{2p}) \times \LGrass(\C^{2q}) \\
\hline
\end{array}
\end{math}
\hfil
}
\end{table}



\def\cftil#1{\ifmmode\setbox7\hbox{$\accent"5E#1$}\else
  \setbox7\hbox{\accent"5E#1}\penalty 10000\relax\fi\raise 1\ht7
  \hbox{\lower1.15ex\hbox to 1\wd7{\hss\accent"7E\hss}}\penalty 10000
  \hskip-1\wd7\penalty 10000\box7} \def\cprime{$'$}
\providecommand{\bysame}{\leavevmode\hbox to3em{\hrulefill}\thinspace}
\renewcommand{\MR}[1]{}
\providecommand{\MRhref}[2]{%
  \href{http://www.ams.org/mathscinet-getitem?mr=#1}{#2}
}
\providecommand{\href}[2]{#2}

\end{document}